\magnification=1200
\input amstex
\documentstyle{amsppt}
\NoBlackBoxes

 \pageheight{9.1truein}
\def\leftitem#1{\item{\hbox to\parindent{\enspace#1\hfill}}}
\def\Vir{\operatorname{Vir}}
\def\dim{\operatorname{dim}}

\def\a{\alpha}
\def\b{\beta}

\def\N{\Bbb N}
\TagsOnRight
\def\supp{\operatorname{supp}}

\def\max{\operatorname{max}}

\def\bC{\Bbb C}

\def\Z{\Bbb Z}
\def\bQ{\Bbb Q}

\def\proclaim#1{\par\vskip.25cm\noindent{\bf#1. }\begingroup\it}
\def\endproclaim{\endgroup\par\vskip.25cm}
\rightheadtext{Classification of weight modules}

\topmatter
\title Classification of irreducible weight
modules over the twisted Heisenberg-Virasoro algebra\endtitle
\author    Rencai Lu, and Kaiming Zhao
\endauthor
\affil
 Institute of Mathematics\\
Academy of Mathematics and System Sciences\\ Chinese Academy of
Sciences\\ Beijing 100080, P.R. China
       \\and\\
       Department of Mathematics\\ Wilfrid Laurier University\\
Waterloo, Ontario\\ Canada N2L 3C5
        Email:  kzhao\@wlu.ca\\
\endaffil
\thanks Research supported by NSERC, and the NSF  of China  (Grants 10371120  and 10431040).
\endthanks

\abstract  In the present paper,
  all irreducible weight modules with finite
dimensional weight spaces over the twisted Heisenberg-Virasoro
algebra are determined. There are two different classes of them.
One class is formed by  simple modules of intermediate series,
whose weight spaces are all $1$-dimensional; the other class
consists of the irreducible highest  weight modules and lowest
weight modules.


\endabstract

\endtopmatter
\document
\smallskip\bigskip

\smallskip\bigskip
\subhead 1. \ \ Introduction
\endsubhead
\medskip

 Let $\bC$ be the  field of complex numbers. The
twisted Heisenberg-Virasoro algebra ${\Bbb L}$ is the universal
central extension of the Lie algebra $\{f(t)\frac{d}{dt}+g(t)| f,g
\in \bC[t,t^{-1}]\}$ of differential operators of order at most
one. More precisely, we have the following definition.

\proclaim{Definition 1.1} The {\bf twisted Heisenbeg-Virasoro
algebra} ${\Bbb L}$ is a Lie algebra  over $\bC$ with the basis
$$\{x_n,I(n),C_D, C_{DI},C_I | n \in \Z\}$$
and the Lie bracket given by
$$[x_n,x_m]=(m-n)x_{n+m}+\delta_{n,-m}\frac{n^3-n}{12}C_D,\tag 1.1$$
$$[x_n,I(m)]=mI(n+m)+\delta_{n,-m}(n^2+n)C_{DI}, \tag 1.2$$
$$[I(n),I(m)]=n\delta_{n,-m}C_I, \tag 1.3$$
$$[\Bbb{L},C_D]=[\Bbb{L},C_{DI}]=[\Bbb{L},C_I]=0. \tag 1.4$$
\endproclaim

Under the projection  $x_n\rightarrowtail t^{n+1}\frac{d}{dt},
I(n)\rightarrowtail t^n$ we can easily see that ${\Bbb L}$ is the
universal central extension of the Lie algebra of differential
operators of order at most one (see [A] for details).

 The highest weight
representations for $\Bbb{L}$ was introduced and studied in [A],
where they precisely  determined the  determinant formula of the
Shapovalov form for Verma modules. In a recent paper [B], Billig
obtained the character formula for irreducible highest weight
modules with trivial action of $C_I$.

The present paper is devoted to determining all irreducible weight
modules with finite dimensional weight spaces over  $\Bbb{L}$.
More precisely we prove that there are two different classes of
them. One class is formed by simple modules of intermediate
series, whose weight spaces are all $1$-dimensional; the other
class consists of the highest(or lowest) weight modules.

The paper is arranged as follows.

In Section 2, we recall some notations and collect known facts
about irreducible, indecomposable modules over the classical
Virasoro algebra.

In Section 3, we determine  all irreducible weight modules of
intermediate series over $\Bbb L$, i.e., irreducible weight
modules with  all $1$-dimensional weight spaces.

In Section 4, determine all irreducible uniformly bounded weight
modules over ${\Bbb L}$ which turn out to be modules of
intermediate series.

In section 5 we obtain the main result of this paper: the
classification of irreducible weight ${\Bbb L}$-modules with
finite dimensional weight space. As we mentioned, they are
irreducible highest, lowest weight modules, or irreducible modules
of the intermediate series.

\subhead 2. \ \ Basics
\endsubhead
\medskip

In this section, we collect some known facts for later use.

The center of ${\Bbb L}$ is four-dimensional and is spanned by
$\{I(0),C_D,C_{DI},C_{I}\}$.

\par

\par For any $e\in\bC$, it is clear that
$$[x_n+eI(n), x_m+eI(m)]=(m-n)(x_{n+m}+eI(n+m)),\,\, \forall n \neq -m,$$
$$[x_n+eI(n),
D(-n)+eI(-n)]=-2n(x_0+eI(0)-eC_{DI}-\frac{e^2}{2}C_I)+\frac{n^3-n}{12}C_D.$$
So $\{x_n+eI(n),x_0+eI(0)-eC_{DI}-\frac{e^2}{2}C_I, C_D| n \in
\Z\backslash \{0\} \}$ span a subalgebra $\Vir[e]$ which is
isomorphic to the classical Virasoro algebra. In many cases, we
shall simply write $\Vir[0]$ as $\Vir$. In this paper we may
sometimes replace $\Vir[0]$ with $\Vir[e]$ by changing the basis
for ${\Bbb L}$.

\par
 Introduce a $\Z$-grading on ${\Bbb L}$ by defining the degrees: deg $x_n$=deg $I(n)$=$n$ and deg
$C_{D}$=deg $C_{DI}$=deg $C_{I}$=0. Set
 \vskip 5pt
$$\Bbb{L}_+=\sum_{n\ge 0}(\bC x_n+\bC I(n)),\,\,\Bbb{L}_-=\sum_{n\le 0}(\bC x_n+\bC I(n)),$$
and
$$\Bbb{L}_{0}=\bC x_0+\bC I(0)+\bC C_D+\bC C_{DI}+\bC C_{I}.$$

\vskip 5pt

A nonzero ${\Bbb L}$-module $V$ is called {\bf trivial} if
$\Bbb{L}V=0$.

For any ${\Bbb L}$-module $V$ and
$(\lambda,\lambda_I,c_D,c_{DI},c_{I})\in \bC^5$, set
$$V_{(\lambda,\lambda_I,c_D,c_{DI},c_{I})}:=\hskip 9cm$$
$$\bigl\{v\in V\bigm|x_0v=\lambda v, I(0)v=\lambda_I,C_D v= c_D v,
C_{DI} v=c_{DI}v,C_I v= c_I v\bigr\},$$ which we generally call
the {\bf weight space} of $V$ corresponding the weight\break
$(\lambda,\lambda_I,c_D,c_{DI},c_{I})$. When
$I(0),C_D,C_{DI},C_{I}$ act as  scalars
$\lambda_I,c_D,c_{DI},c_{I}$ on the whole space $V$, respectively,
we shall simply write $V_{\lambda}$ instead of
$V_{(\lambda,\lambda_I,c_D,c_{DI},c_{I})}$.

\par
An ${\Bbb L}$-module $V$ is called a {\bf weight module} if $V$ is
the sum of all its weight spaces. For a weight module $V$ we
define $$\supp(V):=\bigl\{(\lambda,\lambda_I,c_D,c_{DI},c_{I})\in
\bC^5 \bigm|V_{(\lambda,\lambda_I,c_D,c_{DI},c_{I})}\neq
0\bigr\},$$ which is generally called the {\bf weight set} (or the
{\bf support}) of $V$.

\par
 A nontrivial  weight ${\Bbb L}$-module V is called a {\bf weight module of intermediate
 series} if V is indecomposable and any weight spaces of V is one dimensional.
\par
A  weight ${\Bbb L}$-module V is called a {\bf highest} (resp.
{\bf lowest) weight module} with {highest weight} (resp. {highest
weight}) $(\lambda,\lambda_I,c_D,c_{DI},c_{I})\in \bC^5$, if there
exists a nonzero weight vector $v \in
V_{(\lambda,\lambda_I,c_D,c_{DI},c_{I})}$ such that

 \vskip 5pt 1) $V$ is generated by $v$ as  ${\Bbb L}$-module;

2) $\Bbb{L}_+ v=0 $ (resp. $\Bbb{L}_- v=0 $).
 \vskip 5pt

Let $U:=U(\Bbb L)$ be the universal enveloping algebra of  ${\Bbb
L}$. For any $\lambda,\lambda_I,c_D,c_{DI},c_{I}$ $\in \bC$, let
$I(\lambda,\lambda_I,c_D,c_{DI},c_{I})$ be the left ideal of $U$
generated by the elements $$ \bigl\{x_i,I(i)\bigm|i\in \N
\bigr\}\bigcup\bigl\{x_0-\lambda \cdot 1, I(0)-\lambda_I \cdot 1,
C_D-c_D\cdot 1, C_{DI}-c_{DI}\cdot 1,C_{I}=c_I\cdot 1\bigr\}. $$
Then the {\bf Verma module} with the highest weight
$(\lambda,\lambda_I,c_D,c_{DI},c_{I})$ over ${\Bbb L}$ is defined
as
$$M(\lambda,\lambda_I,c_D,c_{DI},c_{I}):=U/I(\lambda,\lambda_I,c_D,c_{DI},c_{I}).$$
It is clear that $M(\lambda,\lambda_I,c_D,c_{DI},c_{I})$ is a
highest weight module over ${\Bbb L}$ and   contains a unique
maximal submodule. Let $V(\lambda,\lambda_I,c_D,c_{DI},c_{I})$ be
the unique irreducible quotient of
$M(\lambda,\lambda_I,c_D,c_{DI},c_{I})$.

\vskip 5pt
 In the rest of this section, we recall some known facts
about weight representations of the classical Virasoro algebra
which can be considered as a subalgebra of $\Bbb L$:
$$\Vir:= \text{span} {\{x_n,C_D,  | n \in \Z\}}.$$
For details, we refer the readers to [M], the book [KR] and the
references therein.

\par
It is well known that a module of the intermediate series over
$\Vir$ must be one of $V(\a,\b), A(a), B(a)$ for some $ \a,\b,a
\in \bC,$ or one of their quotient submodules, where $V(\a,\b)$
(resp. $A(a),B(a)$) all have basis $\{v_{\a+k}|k\in \Z\}$ (resp.
$\{v_{k}|k\in \Z\}$) such that $C_D$ acts trivially and \vskip 5pt
 $$\hskip -2.5cm V(\a,\b):\,\,x_i v_{\a+k}=(\a+k+\b i)v_{\a+i+k},$$
$$A(a): \,\,x_iv_{k}=(i+k)v_{i+k}, \text{for } k \neq 0;
\,\,x_iv_{0}=i(i+a)v_{i},$$ $$B(a): \,\,x_i v_{k}=kv_{i+k},
\text{for } k \neq -i;\,\, x_iv_{-i}=-i(i+a)v_{0}.$$
 (These facts
appear in many references, for example in [SZ]). We shall
 use $T$ to denote the $1$-dimensional trivial module,
 use $V'(0,0)$ to denote the unique proper nontrivial
submodule of $V(0,1)$ (which is irreducible). \vskip 5pt

Let us recall the extension of one module by another. An
indecomposable module $V$ over $\Vir$ is said to be an {\it
extension of the $\Vir$-module $W_1$ by the $\Vir$-module $W_2$}
if $V$ has a submodule isomorphic to $W_1$ and $V/W_1\simeq W_2$.

\proclaim{Theorem 2.2([MP3])} Let $Z$ be an indecomposable weight
$\Vir$-module  with  weight spaces of dimension  less than or
equal to two.  Then $Z$ is  one of the following:

1) an intermediate series module;

2) an extension of $V(\a,\b)$, $A(a)$ or $B(a)$ by themselves;

3) an extension of $V(\a,\b_1)$ by $V(\a,\b_2)$, where
$\b_1-\b_2=2,3,4,5,6$;

4) an extension of $V(0, \b)$ by $W$, where $\b=2,3,4,5$ and $W$
is one of: $V(0,0)$, $V(0,1)$, $V'(0,0)$, $V'(0,0)\oplus T$, $T$,
$A(a)$ or $B(a)$;

5) an extension of  $T$ by   $V(0,0)$ or $A(a)$, ;

6) an extension of $V'(0,0)$ by  $V(0,0)$ or $A(1)$;

7) an extension  $V(0,1)$ by $V(0,0)$ or $B(0)$;

8) an extension of $V(0,0)$ by   $V(0,1)$ or $A(a)$;

9) an extension of $A(0)$ by  $V(0,0)$;

10) an extension of $B(a)$ by $A(0,1)$;

11) the contragredient extensions of the previous ones;

\noindent where $\a,\b,\b_1,\b_2,a\in\bC$.
\endproclaim

\par
Remark that in the above list, there are some repetitions, and not
all of them can occur.

\proclaim{Theorem 2.3([MP3])}There are exactly two indecomposable
extensions $V=\text{span}\{v_{\a+i},v'_{\a+i}|i\in\Z\}$ of
$V(\a,0)$ by $V(\a,0)$ ($\a \notin \Z$) given by the actions

$$ x_iv_{\a+n}=(\a+n)v_{\a+n+i}, \forall i\in\Z,\,\,\,\text{and} $$

 \noindent  a)
$x_iv_{\a+n}'=(\a+n)v_{\a+n+1}'-iv_{\a+n+1},$ for all $i,n\in\Z$;
or
%

 \vskip 5pt
\noindent  b)  $x_1v_{\a+n}'=(\a+n)v_{\a+n+1}', \hskip 10
ptx_{-1}v_{\a+n}'=(\a+n)v_{\a+n-1}',$ \vskip 5pt
 $
 x_2v_{\a+n}'=(\a+n)v_{\a+n+2}'+\frac{1}{(\a+n+2)(\a+n+1)}v_{\a+n+2},$
\vskip 5pt $
x_{-2}v_{\a+n}'=(\a+n)v_{\a+n-2}'-\frac{1}{(\a+n-2)(\a+n-1)}v_{\a+n-2},$
\vskip 5pt \noindent where $\{v_{\a+n},v'_{\a+n}\}$ forms a base
of $V_{\a+n}$ for all $n\in \Z$.
\endproclaim
\vskip 2pt

 We also need the following result from [MP1]
 \vskip 2pt

\proclaim{Theorem 2.4([MP1, Proposition 3.3]} Suppose that V is a
weight Vir-module with finite dimensional weight spaces. Let $M^+$
(resp. $M^-$ ) be the maximal submodule of V with upper (resp.
lower) bounded weights. If V and V* (the contragredient module of
V) do not contain trivial submodules, there exists a unique
bounded submodule $B$ such that $V=B\bigoplus M^+ \bigoplus M^-$.
\endproclaim

\vskip 5pt \subhead 3. \ \ Irreducible weight modules with weight
multiplicity one
\endsubhead
\medskip

In this section we determine all irreducible weight modules  over
$\Bbb L$ with weight multiplicity one.
\par
\vskip 2pt

Let us first  define a class of ${\Bbb L}$-modules $V(\a,\b; F)$
for $\a,\b,F\in \bC$ as following:    $V(\a,\b;F)$ has base
$\{v_{\a+n}| n\in \Z\}$ and the actions defined by

$$x_iv_{\a+k}=(\a+k+\b i)v_{\a+i+k}, \,\,\,I(i)v_{\a+k}=Fv_{\a+k+i} ,$$
$$C_Dv_{\a+k}=0 , \,\,\, C_Iv_{\a+k}=0, \,\,\,  C_{DI}v_{\a+k}=0\,\,\forall\,\, i,k \in
\Z.$$

\vskip 5pt

It is easy to check that the ${\Bbb L}$-module $V(\a,\b;F)$ is
reducible if and only if $F=0$, $\a \in \Z$ and $\b=0,1$. Denote
the unique (isomorphic) nontrivial sub-quotient module of
$V(0,0;0)$ and $V(0,1;0)$ as $V'(0,0;0)$.

\vskip 5pt

Let V be a nontrivial irreducible weight $\Bbb L$-module with
weight multiplicity one. Then we may assume that
 $I(0), C_{DI},  C_{I}, C_{D}$ act as scalars $F,  c_{DI}, c_{I}, c_{D}$
 respectively. It is clear that $c_D=0$ and $\supp(V) \subset \a+\Z$ for some $\a\in\bC$.

 \proclaim{Lemma 3.1} If $V\simeq V(\a,\b)$ as $\Vir[0]$-module,
 then $V\simeq V(\a,\b;F)$ as ${\Bbb L}$-module for some $F\in \bC$.
\endproclaim
\demo{Proof} Suppose that $V=\oplus_{k\in\Z}\bC v_k$ and the
action is given by $x_iv_n=(\a+n+i\b)v_{n+i}$. Let $I(i)
v_{\a+n}=F_{i,n}v_{\a+n+i}$ where $F_{i,n}$ are constants.
Computing $[x_i,I(j)]v_{n+\a}$, and $[I(i),I(j)]v_{n+\a}$ gives

 $$F_{0,n}=F,\eqno(3.1)$$
 $$F_{j,n+i}F_{i,n}-F_{i,n+j}F_{j,n}=j\delta_{i,-j}c_I,\eqno(3.2)$$
$$(\a+n+j+i\b)F_{j,n}-F_{j,n+i}(\a+n+i\b)=jF_{i+j,n}+\delta_{i,-j}(i^2+i)c_{DI},\eqno(3.3)$$
for all $i,j,n\in\Z$.
\medskip
With $j=1$, (3.3) becomes (even for $i=1$)
$$F_{i+1,n}=(\a+i\b+n+1)F_{1,n}-F_{1,n+i}(\a+i\b+n),$$
i.e.,
$$F_{i,n}=(\a+(i-1)\b+n+1)F_{1,n}-(\a+n+(i-1)\b)F_{1,n+i-1}.\eqno(3.4)$$

Taking $i=0$ in (3.4), we have
 $F_{0,n}=(\a-\b+n+1)F_{1,n}-(\a-\b+n)F_{1,n-1},$ i.e.,
$$(\a-\b+n+2)F_{1,n+1}=(\a-\b+n+1)F_{1,n}+F.\eqno(3.5)$$

 \vskip 2pt
With $i=2$, (3.4) gives
$$F_{2,n}=(\a+\b+n+1)F_{1,n}-(\a+n+\b)F_{1,n+1}.\eqno(3.6)$$  Combining the above two formulas we deduce that
\vskip 5pt

$(\a-\b+n+2)F_{2,n}$

$=(\a-\b+n+2)(\a+\b+n+1)F_{1,n}-(\a+n+\b)((\a-\b+n+1)F_{1,n}+F)$

$=[(\a-\b+n+2)(\a+\b+n+1)-(\a+n+\b)(\a-\b+n+1)]F_{1,n}-(\a+n+\b)F,$

i.e.,
 $$(\a-\b+n+2)F_{2,n}=(2\a+2n+2)F_{1,n}-(\a+n+\b)F.\eqno(3.7)$$

\vskip 5pt If $i=1, j=2$ in (3.2),  we have
$F_{2,n+1}F_{1,n}-F_{1,n+2}F_{2,n}=0$, i.e.,
$$(\a-\b+n+3)F_{2,n+1}(\a-\b+n+2)F_{1,n}=(\a-\b+n+3)F_{1,n+2}(\a-\b+n+2)F_{2,n}.\eqno(3.8)$$

Applying (3.7) and (3.5) to (3.8) we deduce

$$[(2\a+2n+4)((\a-\b+n+1)F_{1,n}+F)-(\a-\b+n+2)(\a+\b+n+1)F]F_{1,n}$$
$$=[(\a-\b+n+1)F_{1,n}+2F][(2\a+2n+2)F_{1,n}-(\a+n+\b)F],$$
simplifying to give
 $$(F_{1,n}-F)((\a-\b+n+1){F_{1,n}}-(\a+n+\b) F)=0.\eqno(3.9)$$
\medskip

{\it Case 1: $F=0$.}
\medskip
It is clear from (3.9) that $F_{1,n}=0$ for all $n \neq \b-\a-1$.
Next we want to show that $F_{1,n}=0$ for all $n$.

\par Suppose  that $\b-\a \in \Z$ and
$F_{1,\b-\a+1} \neq 0$. We may assume that $\a=\b$, and $F_{1,-1}
\neq 0$. From (3.4) we see that
 $$F_{j,n}=0,\,\,\,\, \forall (j,n) \notin \{ (i,-1), (i,-i) | i
 \in \Z\};$$
 $$F_{j,-1}=j\b F_{1,-1},\,\,\,\, \forall j \neq 1;$$
 $$F_{n,-n}=-(n\b-n)F_{1,-1}, \forall n \neq 1.$$

 If $n=-j=1$ and $i=2$ in (3.3),  using the above formulas
 we deduce $(\b-1)\b=0$.
 So we have $\a=\b=0$ or $\a=\b=1$.

 In the case  $\a=\b=0$, noting that
  $F_{j,0}=0$ for all $j$, we know that $V_{0}$ is an  ${\Bbb L}$ submodue of V, a contradiction.
Thus $F_{1,n}=0$ for all $n$.

 If $\a=\b=1$,  (noting that $V(\a,\b)^*\simeq V(1-\a, 1-\b)$) the contragredient module of $V$ contains a
 $1$-dimensional trivial ${\Bbb L}$-submodule, again a contradiction.
Thus $F_{1,n}=0$ for all $n$ either.

 In both cases,
using (3.4) we deduce that $F_{m,n}=0$ for all $m, n\in\Z$. From
(3.2) and (3.3) we can easily deduce that $c_{I}=c_{DI}=0$. Hence
$V\simeq V(\a,\b;0)$ in Case 1.

\medskip
{\it Case 2: $F\neq 0$.}
\medskip

{\it Case 2.1: $F_{1,n}=F$ for some $n$, say $F_{1,0}=F$.}

 From (3.5) we deduce that
$(\a-\b+n+1+k)F_{1,n+k}=(\a-\b+n+1)F_{1,n}+kF$ for all $k$, then
 $F_{1,n}=F$ if $\a-\b+n+1\ne0$.

If $F_{1,n_0}\ne F$ for $n_0=\b-\a-1\in\Z$, from (3.9) we see that
$\a+\b+n_0=0$, i.e., $\b=1/2$ and $\a=-n_0-1/2$.

From (3.4) we see that
 $$F_{j,n}=F,\,\,\,\, \forall (j,n) \notin \{ (i,n_0), (i,n_0+1-i) | i
 \in \Z\};$$
 $$F_{j,n_0}=(j/2)F_{1,n_0}+(1-j/2) F,\,\,\,\, \forall j \neq 1.\eqno(3.10)$$
With $n=n_0$, $i\ne j$, $ij\ne0$ and $i+j\ne1$, (3.2) yields
$F_{i,n_0}=F_{j,n_0}$. Then using (3.10) we deduce that
$F_{1,n_0}=F$. So we have $F_{j,n_0}=F$ for all $j\in\Z$. Applying
(3.4) again we see that $F_{m,n}=F$ for all $m, n\in\Z$. From
(3.2) and (3.3) we can easily deduce that $c_{I}=c_{DI}=0$. Hence
$V\simeq V(\a,\b;F)$ in this case.

\vskip 5pt {\it Case 2.2:  $F_{1,n}\ne F$ for all $n\in\Z$.}

 Then from (3.9) we see that
 $$(\a-\b+n+1){F_{1,n}}=(\a+n+\b) F,\,\,\forall\,\,\,n
\in \Z.\eqno(3.11)$$
\medskip
Combining this with (3.4) we obtain that
$$\hskip -6cm (\a-\b+n+1)(\a-\b+n+j)F_{j,n}=$$
$$\hskip -3cm (\a+(j-1)\b+n+1)(\a+\b+n)(\a-\b+n+j)F$$
$$-(\a+n+(j-1)\b)(\a+\b+n+j-1)(\a-\b+n+1)F,\,\,\forall \,\,\,n\in\Z.\eqno(3.12)$$

Multiplying (3.3) by
$(\a-\b+n+j)(\a-\b+n+1)(\a-\b+n+i+1)(\a-\b+n+i+j)$,  using (3.12)
we obtain that
$$\hskip -3cm (\a-\b+n+i+1)(\a-\b+n+i+j)(\a+n+j+i\b)$$
$$\cdot (\a+(j-1)\b+n+1)(\a+\b+n)(\a-\b+n+j)
$$
$$\hskip -3cm -
  (\a-\b+n+i+1)(\a-\b+n+i+j)(\a+n+j+i\b )$$
  $$\hskip 2cm\cdot (\a+n+(j-1)\b)(\a+\b+n+j-1)(\a-\b+n+1)$$
  $$\hskip -4cm -
  (\a-\b+n+j)(\a-\b+n+1)(\a+n+i\b)$$
  $$\hskip 2cm\cdot (\a+(j-1)\b+n+i+1)(\a+\b+n+i)(\a-\b+n+i+j)$$
  $$\hskip -4cm +
  (\a-\b+n+j)(\a-\b+n+1)(\a+n+i\b)$$
  $$\hskip 2cm\cdot (\a+n+i+(j-1)\b)(\a+\b+n+i+j-1)(\a-\b+n+i+1)$$
  $$-
  j(\a-\b+n+j)(\a-\b+n+i+1)(\a+(i+j-1)\b+n+1)(\a+\b+n)(\a-\b+n+i+j)$$
  $$\hskip -2.5cm +
  j(\a-\b+n+j)(\a-\b+n+i+1)(\a+n+(i+j-1)\b)$$ $$\hskip 3cm\cdot (\a+\b+n+i+j-1)(\a-\b+n+1)=0$$
 for all $i+j \neq 0$, simplifying to give
$$ij\b \left( 2\,\b-1 \right)  \left( \b-1 \right)  \left( j-1 \right)
 \left( 1+i \right)  \left( -2\,\b+i+j+2\,\a+2\,n+1 \right)=0,\,\,\,\forall i,j,n\in\Z.$$
Hence  $\b=0$ , $\b=\frac{1}{2}$, or $\b=1$. Next we discuss them
separately.
\medskip
{\it Case 2.2.1:   $\b=0$.}
\medskip
If $\a \in \Z$, we may suppose  that $\a=0$. From (3.12) with
$n=0$  we deduce that $F_{j,0}=F$ for all $j\ne0$. So  that V
contains a nonzero trivial ${\Bbb L}$-submodule $V_{0}$. This
contradiction tells that this case do not occur.

Now suppose  that $\a \notin \Z$, from (3.12) we have
 $F_{j,n}=\frac{\a+n}{\a+n+j}F.$
 Let $w_{\a+n}=\frac{v_{\a+n}}{\a+n}.$ Then
 $$x_iw_{\a+n}=\frac{(\a+n)}{\a+n}v_{\a+n+i}=(\a+n+i)w_{\a+n+i},$$
$$I(i)w_{\a+n}=\frac{Fv_{\a+n+i}}{\a+n+i}=Fw_{\a+n+i}.$$
From (3.2) and (3.3) we can easily deduce that $c_{I}=c_{DI}=0$.
Thus $V\simeq V(\a,1;F)$ as ${\Bbb L}$-module.
\medskip
{\it  Case 2.2.2:   $\b=1$}
\medskip
If $\a \in \Z$, we may assume that $\a=0$. From (3.11) we have
$nF_{1,n}=(n+1)F$ for all $n \in \Z$, in particularly, $F=0$, a
contradiction.

Now suppose  that $\a \notin \Z.$ we have
 $$F_{j,n}=(\a+j+n)\frac{\a+1+n}{\a+n}c-(\a+n+j-1)\frac{\a+n+j}{\a-1+n+j}F=\frac{\a+n+j}{\a+n}F.$$
 Let $w_{\a+n}=(\a+n)v_{\a+n}$. By a similar discussion as in Case 2.2.1, $V\simeq V(\a,0;F)$.

\medskip
{\it  Case 2.2.3:  $\b=\frac{1}{2}$.}
\medskip
If $\a \notin \Z +1/2$, from (3.11) we know that  $F_{1,n}=F$ for
all $n$, contradicting the assumption. So this case does not
occur.

If $\a \in \Z +1/2$, we may assume that $\a=1/2$.
 Hence we have $F_{1,n}=F$ for all $n \neq -1$,
contradicting the assumption. So this case does not occur. This
completes the proof. \hfill $ $\qed
\enddemo

\proclaim{Lemma 3.2} If $V \ncong V(\a,\b)$ as $Vir[0]$-module,
then $V \simeq V'(0,0;0)$.

\endproclaim

\demo{Proof} Recall that $\supp(V)\subset \a+\Z$  and $c_D=0$.
Since $V \ncong V(\a,\b)$ as $\Vir[0]$-module, then $V$ as
$\Vir[0]$-module (also as $\Vir[e]$-module for any $e\in\bC$) is
possibly isomorphic to: $V'(0,0), V'(0,0)\oplus T,$ $A(a)$ or
$B(a)$. Then the action of the element
$D(0)+eI(0)-eC_{DI}-(e^2/2)C_I$ has integer eigenvalues, i.e.,
$eF-ec_{DI}-(e^2/2)c_I$ are integers for all $e$. Thus we have
$F=c_{DI}$, and $c_I=0$.

We may assume that $V=\oplus_{k\in\Z}\bC v_k$ (only $v_0$ can be
zero) and the action is given by
$$x_iv_n=nv_{n+i}, \,\,\,\,\,I(i) v_{n}=F_{i,n}v_{n+i},\,\,\,\,\,{\text{ for} }  \,\,\,n(n+i)\ne0,$$
where $F_{i,n}$ are some constants.

The actions $[I(i),I(j)]v_n$ and $[x_i,I(j)]v_n$ give $$F_{0,n}=F,
\,\,\,\,\,{\text{ for} } \,\,\, n \neq 0,\eqno(3.13)$$
$$F_{j,n+i}F_{i,n}-F_{i,n+j}F_{j,n}=0,\,\,\,\,\,{\text{ for} }  \,\,\,
 n(n+i+j) (n+i)(n+j) \neq 0,\eqno(3.14)$$
$$(n+j)F_{j,n}-nF_{j,n+i}=jF_{i+j,n}+\delta_{i,-j}(i^2+i)F,\,\,\,\,\,{\text{ for} }
 \,\,\,n(n+i+j)(n+i)(n+j)\neq 0.\eqno(3.15)$$ \vskip 5pt

Letting $i=1$ in (3.14), we obtain that
$$F_{j,n+1}F_{1,n}=F_{1,n+j}F_{j,n}
,\,\,\,\,\,{\text{ for} }  \,\,\,
 n(n+1+j) (n+1)(n+j) \neq 0.\eqno(3.16)$$
 Letting $j=1$ in (3.15), we obtain that
 $$F_{i+1,n}=(n+1)F_{1,n}-nF_{1,n+i}-\delta_{i,-1}(i^2+i)F,\,\,\,\,\,{\text{ for} }
 \,\,\,n(n+i+1)(n+i)(n+1)\neq 0,$$
i.e., $$F_{j,n}=(n+1)F_{1,n}-nF_{1,n+j-1},\,\,\,\,\,{\text{ for} }
 \,\,\,n(n+1)(n+j)(n+j-1)\neq 0.\eqno(3.17)$$
 Letting $j=0$, we obtain that
$$F_{0,n}=(n+1)F_{1,n}-nF_{1,n-1},\,\,\,\,\,{\text{ for} }
 \,\,\,n(n+1)(n-1)\neq 0,$$
i.e.,
 $$(n+2)F_{1,n+1}=(n+1)F_{1,n}+F\,\,\,\,\,{\text{ for} }
 \,\,\,n(n+1)(n+2)\neq 0.\eqno(3.18)$$
Thus we have
$$(n+3)F_{1,n+2}=(n+1)F_{1,n}+2F,\,\,\,\,\,{\text{ for} }
 \,\,\,n(n+1)(n+2)(n+3)\neq 0.$$
Letting $j=2$ in (3.17), we have
$$F_{2,n}=(n+1)F_{1,n}-nF_{1,n+1}, \,\,\,\,\,{\text{ for} }
 \,\,\,n(n+1)(n+2)\neq 0.\eqno(3.19)$$
Using this we deduce  that
$$(n+2)F_{2,n}=(2n+2)F_{1,n}-nF,\,\,\,\,\,{\text{ for} }
 \,\,\,n(n+1)(n+2)(n+3)\neq 0.\eqno(3.20)$$
With $j=2$ in (3.16), we know that $$(n+2)(n+3)F_{2,n+1}F_{1,n}=
(n+2)F_{2,n}(n+3)F_{1,n+2}.$$ Applying (3.18) and (3.20) we deduce
that $$[(2n+4)((n+1)F_{1,n}+F)-(n+2)(n+1)F]F_{1,n}=$$
$$[(n+1)F_{1,n}+2F][(2n+2)F_{1,n}-nF],\,\,\,\,\,{\text{ for} }
 \,\,\,n(n+1)(n+2)(n+3)\neq 0,$$
simplifying to give
 $$(F_{1,n}-F)((n+1){F_{1,n}}-n F)=0,\,\,\,\,\,{\text{ for} }
 \,\,\,n(n+1)(n+2)(n+3)\neq 0.$$
Then, for any $n$ with $n(n+1)(n+2)(n+3)\neq 0$, $F_{1,n}=F$  or
 $(n+1){F_{1,n}}=n F.$
Using (3.18) we deduce that
$$F_{1,n}=F, \,\,\, \forall \,\,\, n \geq 1\,\,\,\text{or}\,\,\,
F_{1,n}=\frac{n}{n+1}F,\,\,\, \forall \,\,\, n \geq
1,\eqno(3.21)$$ and
$$F_{1,n}=F, \,\,\, \forall \,\,\, n \leq
-4\,\,\,\text{or}\,\,\,F_{1,n}=\frac{n}{n+1}F, \,\,\, \forall
\,\,\, n \leq -4.\eqno(3.22)$$

Applying to (3.19), we have $F_{2,n}=F$ for all $n \geq 1$ or
$F_{2,n}=\frac{n}{n+2}F$ for all $n \geq 1$. \vskip 5pt

Take $i=-2, j=2$ and $n>4$ in (3.15), we get $F=0$. Then (3.21)
and (3.22) infer that $$F_{1,n}=0,\,\,\text{if}\,\,n>0
\,\,\,\text{or}\,\,\, n<-3.$$ Applying this to (3.18) with $n=-4$
and $-3$, we obtain that
$$F_{1,n}=0,\,\,\text{if}\,\,n \neq 0,-1.\eqno(3.23)$$
\medskip
Note that ${\Bbb L}$ is generated by $\{X_n, I(1)| n\in \Z\}$ as
Lie algebra. If  $V\simeq V'(0,0)$ as $\Vir$-module, then $V\simeq
V'(0,0;0)$.

Next we assume that $v_0\ne0$.

Applying (3.23) to (3.17) we obtain that

$$I(k)v_n=0,\,\,\,\,\,{\text{ for} }
 \,\,\,n(n+1)(n+k)(n+k-1)\neq 0.\eqno(3.24)$$
Noting that $I(0)v_{n}=0$,  we have
$$\hskip -3cm I(k)v_{-1}=\frac{1}{k-2}I(k)x(-k+1)v_{k-2}$$
$$=\frac{1}{k-2}x(-k+1)I(k)v_{k-2}=0 \,\,{\text{ for}
}\,\,(k-1)(k-2)\neq 0;\eqno(3.25)$$

$$I(k)v_{-k+1}=I(k)x(-k)v_{1}=x(-k)I(k)v_{1}=0,\,\,\,\,\,{\text{ for} }\,\,\, (k+1)(k-1)\neq 0. \eqno(3.26)$$

From (3.24)-(3.26) we have $$I(j)v_{n}=0 \,\,\,\,\,{\text{ for}}
\,\,\,\, n(n+j)\neq 0.\eqno(3.27)$$

 Note that $F=c_{DI}=c_{I}=0$.

 If $V_0$ is a submodule over
 $\Vir[0]$, we write $I(i)v_0=a_iv_i$ where $a_i\in\bC$.
 From the action $[x_k, I(i)]v_0=iI(i+k)v_0$ and using (3.27) and the fact $x_kv_0\in\bC v_k$,
  we deduce that $a_i=a$, a fixed
 constant, for all $i\ne0$. If $a=0$ we see that $V_0$ is an $\Bbb L$-submodule, a
 contradiction. So $a\ne0$.


Replacing $\Vir[0]$ by $\Vir[1]$,
 we may assume  that  $V_0$ is not a $\Vir[0]$-submodule.
So we can always assume that $V_0$ is not a $\Vir[0]$-submodule.
Hence $V\simeq A(a)$ as $\Vir[0]$-module
 for some $a\in\bC$. Computing the action
 $[x_k,I(i)]v_{-i}=iI(i+k)v_{-1}=0$ for $k$ with $k(i+k+1)\ne0$, we deduce that
 $x_kI(i)v_{-i}=0,$ that is
 $I(i)v_{-i}=0$ for all $i$. So $\oplus_{i\ne0}\bC v_i$ is a
 proper submodule over $\Bbb L$, a
 contradiction. This case do not occur. This completes the proof.
\hfill $ $\qed
\enddemo

Combining the previous two lemmas, we obtain

\proclaim{Theorem 3.3} Suppose  that  V is a nontrivial
irreducible weight $\Bbb L$-module with weight multiplicity one.
Then we have $V\simeq V(\a,\b;F)$ or $V\simeq V'(0,0;0)$ for some
$\a,\b, F\in\bC$.
\endproclaim

\subhead 4. \ \ Uniformly bounded irreducible weight modules
\endsubhead
\medskip

In this section,  we assume that $V$ is a uniformly bounded
nontrivial irreducible weight module over $\Bbb L$. So there
exists $\a\in\bC$ such that $\supp(V)\subset \a+\Z$. From
representation theory of $\Vir$, we have $C_D=0$ and $\dim
V_{\a+n}=p$ for all $\a+n \neq 0$. If $\a \in \Z$, we also assume
 that $\a=0$.

Consider $V$ as a $\Vir$-module. We have a $\Vir$-submodule
filtration
    $$0=W^{(0)}\subset W^{(1)} \subset W^{(2)}\subset \cdots \subset W^{(p)}=V,$$
where $W^{(1)}, \cdots ,W^{(p)}$ are $\Vir[0]$-submodules of $V$,
and the quotient modules \break $W^{(i)}/W^{(i+1)}$ have weight
multiplicity one for all nonzero weights.

Choose $v_n^1, \cdots, v_n^p \in V_{\a+n}$ such that the images of
 $v_n^i+W^{(i-1)}$ form a basis of $(W^{(i)}/W^{(i-1)})_{\a+n}$ for all $\a+n\neq 0$. We may suppose that
 $$x_i
 (v_n^1, \cdots ,v_n^p)=
 (x_iv_n^1, \cdots ,x_iv_n^p)=(v_{n+i}^1, \cdots ,v_{n+i}^p)A_{i,n},$$ where $A_{i,n}$
 are upper triangular $p\times p$ matrices, and $A_{i,n}(j,j)=\a+n+i\b_j $.
 Denote $$I(i)
 (v_n^1, \cdots ,v_n^p)=(v_{n+i}^1, \cdots ,v_{n+i}^p)F_{i,n},$$
where $F_{i,n}$
 are $p\times p$ matrices.

 The Lie brackets gives
\vskip -.3cm
$$\hskip -9.5cm F_{0,n}=FI_p, \tag 4.1 $$

\vskip -.5cm

$$\hskip -5.8cm F_{i,j+n}F_{j,n}-F_{j,i+n}F_{i,n}=j\delta_{i,-j}c_II_p,\tag 4.2$$
\vskip -.5cm

$$\hskip -5.1cm A_{i,j+n}A_{j,n}-A_{j,i+n}A_{i,n}=(j-i)A_{i+j,n},\tag 4.3 $$
\vskip -.5cm

$$\hskip -3.1cm A_{i,j+n}F_{j,n}-F_{j,i+n}A_{i,n}=jF_{i+j,n}+\delta_{i,-j}(i^2+i)c_{DI}I_p,\tag 4.4 $$
where $I_p$ is the identity matrix, and the last three formulas
have the restriction $(\a+n)(\a+n+i)(\a+n+j)(\a+n+i+j)\ne0$. We
shall denote the $(i,j)$-entry of a matrix $A$ by $A(i,j)$.

\medskip
\proclaim{Lemma 4.1} If all nontrivial irreducible sub-quotient
$\Vir[0]$-modules of V are isomorphic to $V'(0,0)$ and $F=C_{DI}$,
 then $V \simeq V'(0,0,F)$.
\endproclaim
\demo{Proof} From the assumption and Theorem 2.1 we may choose
$W^{(1)}$ such that $\dim (W^{(1)})_0 \leq 1$ (if $V$ contains a
trivial submodule $\bC v _0$, then span$\{v_1^0,
I(k)v_1^0|k\in\Z\}$ is a Vir-submodule which can be chosen as
$W^{(1)}$).

\medskip
 {\bf Claim}. {\it The $(k,1)$-entry $F_{j,n}(k,1)=0$ for all $k\geq 2$, $n \neq 0$ and
 $j+n \neq 0$.}
\medskip
 {\it Proof of Claim.} Suppose  that we have $F_{j,n}(k,1)=0$ for all $k\geq k_0+1(k_0\geq 2)$, $n \neq 0$ and
 $j+n \neq 0$. We only need prove that $F_{j,n}(k_0,1)=0$ for all $n \neq 0$ and
 $j+n \neq 0$.

The $(k_0,1)$-entry of (4.4) gives
 $$(n+j)F_{j,n}(k_0,1)-F_{j,i+n}(k_0,1)n=jF_{i+j,n}(k_0,1),\,\,\,\text{if}\,\,\, n\ne 0,-i,-j,-i-j.\eqno(4.5)$$
Letting $j=1$ in (4.1), we have the $(k_0,1)$-entry
 $$F_{i+1,n}(k_0,1)=(n+1)F_{1,n}(k_0,1)-nF_{1,i+n}(k_0,1),\,\,\,{\text{if}}\,\,\,n \neq 0,-1,-i,
 -1-i,$$
i.e.,
 $$F_{j,n}({k_0},1)=(n+1)F_{1,n}({k_0},1)-nF_{1,n+j-1}({k_0},1),\,\,\,{\text{if}}\,\,\,n \neq 0,-1,-j+1,
 -j.\eqno(4.6)$$
 Letting $j=0$ in (4.6),  we have
 $0=(n+1)F_{1,n}({k_0},1)-nF_{1,n-1}({k_0},1)$, for all $n \neq
 0,\pm1$,
 i.e.,
$$(n+2)F_{1,n+1}({k_0},1)=(n+1)F_{1,n}({k_0},1),\,\,\,{\text{if}}\,\,\,n \neq
0,-1,-2.$$
 Hence $F_{1,n}({k_0},1)=\frac{2}{n+1}F_{1,1}({k_0},1)$
for all $n\geq 1$, and
$F_{1,n}({k_0},1)=\frac{-1}{n+1}F_{1,-2}({k_0},1)$ for all $n\le
-2$.

 Suppose  that $F_{1,1}({k_0},1)\neq 0$, i.e., $0 \neq I(1)v_1^1+W^{(k_0-1)} \in
W^{(k_0)}/W^{(k_0-1)}$.

 If $W^{(1)}$ contains $V'(0,0)$ as  a Vir-submodule, then we have a nonzero
element  $I(1)v_1^1 +W^{(k_0-1)}$ in $(W^{(k_0)}/W^{(k_0-1)})_2$
such that
$x_{-1}I(1)v_1^1+W^{(k_0-1)}=(I(0)+I(1)x_{-1})v_1^1+W^{(k_0-1)}=W^{(k_0-1)}$,
i.e., the Vir-module $W^{(k_0)}/W^{(k_0-1)}$ has a nontrivial
submodule not isomorphic to $V'(0,0)$, contradicting the
assumption in the lemma.

If $W^{(1)}$ contains a trivial submodule $\bC v_0$, then we can
replace $W^{(1)}$ by \break span$\{v_1^0, I(k)v_1^0|k\in\Z\}$.
Then we have a nonzero element  $I(1)v_1^1 +W^{(k_0-1)}$ in
$(W^{(k_0)}/W^{(k_0-1)})_2$ such that
$x_{-1}I(1)v_1^1+W^{(k_0-1)}=(I(0)+I(1)x_{-1})v_1^1+W^{(k_0-1)}=W^{(k_0-1)}$,
i.e., the Vir-module $W^{(k_0)}/W^{(k_0-1)}$ has a nontrivial
submodule not isomorphic to $V'(0,0)$, contradicting the
assumption in the lemma. Hence

$$F_{1,1}({k_0},1)=0.$$ Similarly we have
$$F_{1,-2}({k_0},1)=0.$$

Applying these to (4.6) we deduce that $F_{i,n}({k_0},1)=0$ for
all $n \neq 0, -1,-i,-i+1$. Letting $i=-n-1$ in (4.5) for suitable
$n$ we deduce that $F_{i,-1}({k_0},1)=0$, and letting $n=-j+1$ in
(4.5) for suitable $i$ we deduce that $F_{i,-i+1}({k_0},1)=0$.
 So we have proved this Claim.

\medskip
This claim ensures that $I(i)v^1_j\in W^{(1)}$ if $j(i+j)\ne0$.
Consider the action of $\Bbb L$ on $W^{(1)}$. By the same argument
as in Lemma 3.2 (from the second paragraph to (3.23) in the proof,
where all $F_{i,j}$ should be replaced by $F_{i,k}({1},1)$) we
have $F=0$. Applying these to (3.23) to (3.17) we deduce that
$F_{i,n}(1,1)=0$ for all $n \neq 0, -1,-i,-i+1$. Letting $i=-n-1$
in (3.15) for suitable $n$ we deduce that $F_{i,-1}(1,1)=0$, and
letting $n=-j+1$ in (3.15) for suitable $i$ we deduce that
$F_{i,-i+1}(1,1)=0$. Hence we obtain that $F_{j,n}(1,1)=0$ for all
$j+n \neq 0$ and $n\neq 0.$ \vskip 5pt

By computing the actions $[x_i,I(j)]v^1_{-i-j}=jI(i+j)v^1_{-i-j}$
we deduce that \break $\dim \sum_{j\in \Z }\bC I(j)v^1_{-j} \leq
1$. It is clear that $I(k)I(j)v^1_{-j}=I(j)I(k)v^1_{-j}=0$ for all
$k\ne j$, that $(j-k)I(j)I(j)v^1_{-j}=[x_k,I(j-k)]I(j)v^1_{-j}=0$
for many suitable $k$. Hence all weight spaces of $W=U({\Bbb
L})v_{1}^p$ is one dimensional. Combining with Lemma 3.3, we have
proved this lemma. \hfill $ $\qed
\enddemo

\proclaim{Lemma 4.2} If any nontrivial irreducible sub-quotient
$\Vir[0]$-module of V is  isomorphic to  $V(\a,0)$, where $\a
\notin \Z$, then we have $V\simeq V(\a,0,F)$ for some $F\in\bC$.
\endproclaim
\demo{Proof} We use the same notations and similar discussions as
in the proof of Lemma 4.1. Suppose  that we have $F_{j,n}(k,1)=0$
for all $k\geq k_0+1(k_0\geq 2)$, $n \neq 0$ and
 $j+n \neq 0$. We first want to prove that $F_{j,n}(k_0,1)=0$ for all $n $ and
 $j$.

The $(k_0,1)$-entry of (4.4) gives
 $$(\a+n+j)F_{j,n}(k_0,1)-F_{j,i+n}(k_0,1)(\a+n)= j F_{i+j,n}(k_0,1).\eqno(4.7)$$
Letting $j=1$ in (4.7), we have the
 $$ F_{i+1,n}(k_0,1)=(\a+n+1)F_{1,n}(k_0,1)-(\a+n)F_{1,i+n}(k_0,1),$$
i.e.,
 $$ F_{j,n}({k_0},1)=(\a+n+1)F_{1,n}({k_0},1)-(\a+n)F_{1,n+j-1}({k_0},1).\eqno(4.8)$$
 Letting $j=0$ in (4.8),  we have
 $0=(\a+n+1)F_{1,n}({k_0},1)-(\a+n)F_{1,n-1}({k_0},1)$.
 Hence $$F_{1,n}({k_0},1)=\frac{\a+ 1}{\a+n+1}F_{1,0}({k_0},1)
 ,\,\,\,\forall\,\,\,n\in\Z.$$
Applying to (4.8) we obtain that
$$F_{j,n}({k_0},1)=\frac{(\a+1)j}{\a+n+j}F_{1,0}(k_0,1),\,\,\,\forall\,\,\,j, n\in\Z.\eqno(4.9)$$
Suppose  that $F_{1,0}(k_0,1)\ne0$. By re-scalaring
$\{v_i^{k_0}|i\in\Z\}$ we may assume that $$
F_{1,0}(k_0,1)=1/(\a+1).\eqno(4.10)$$
\medskip

{\it Case 1: $k_0\geq 3$.}
\medskip
{\it Case 1.1: $W^{(k_0)}/W^{(k_0-2)}$  is decomposable  over
$\Vir$.}
\medskip
In this case we can suitable choose $\{v^k_j|k,j\in\Z\}$ so that
besides (4.9) we also have
$$F_{j,n}({k_0-1},1)=\frac{(\a+1)j}{\a+n+j}F_{1,0}(k_0-1,1),\,\,\,\forall\,\,\,j, n\in\Z.\eqno(4.11)$$
If $F_{1,0}(k_0 ,1)\ne0$, we know that $I(1)v^1_0\mod W^{(k_0-2)}
,v^{k_0}_1\mod W^{(k_0-2)}$ are linearly independent, and that
$I(1)v^1_0\mod W^{(k_0-2)} ,v^{k_0-1}_1\mod W^{(k_0-2)}$ are
linearly independent. Then we can re-choose $W^{(k_0-1)}$ and
$\{v^{k_0-1}_j| j\in\Z\}$ such that $v^{k_0-1}_1=I(1)v^1_0.$ Then
$F_{1,0}({k_0},1)=0$, furthermore
 $F_{j,n}({k_0},1)=0 ,$ for all $j, n\in \Z.$
\medskip
{\it Case 1.2:    $W^{(k_0)}/W^{(k_0-2)}$ is indecomposable  over
$\Vir$.}
\medskip
From Theorem 2.3, we need consider two subcases. \vskip 5pt

{\it Case 1.2.1:  $A_{i,n}(k_0-1,k_0)=-i$ for all $i, n \in \Z$.}

 Using (4.9) and (4.10), from the $(k_0-1,1)$-entry of (4.4), we
 obtain
$$(\a+n+j)F_{j,n}(k_0-1,1)-\frac{ ij}{\a+n+j} -F_{j,i+n}(k_0-1,1)(\a+n) =jF_{i+j,n}(k_0-1,1).\eqno(4.12)$$
 Letting $j=1$ and $i=-1$, we obtain that
$$(\a+n+1)F_{1,n}(k_0-1,1)+\frac{1}{\a+n+1}-F_{1,n-1}(k_0-1,1)(\a+n,)=0,$$
 i.e.,
  $$F_{1,n}(k_0-1,1)=\frac{\a+n}{\a+n+1}F_{1,n-1}(k_0-1,1)-\frac{1}{(\a+n+1)^2}.\eqno(4.13)$$
  Letting $i=j=1$, and $j=-i=2$ in (4.12) respectively, we obtain that
  $$F_{2,n}(k_0-1,1)=(\a+n+1)F_{1,n}(k_0-1,1)-\frac{1}{\a+n+1}-F_{1,1+n}(k_0-1,1)(\a+n),$$
  $$F_{2,n}(k_0-1,1)=\frac{\a+n}{\a+n+2}F_{2,n-2}(k_0-1,1)-\frac{4}{(\a+n+2)^2}.$$
 Using the above two formulas and (4.13) we deduce that
$$F_{2,n+2}(k_0-1,1)=\frac{\a+n+2}{\a+n+4}F_{2,n}(k_0-1,1)-\frac{4}{(\a+n+4)^2}$$
  $$\hskip -4cm =\frac{\a+n+2}{\a+n+4}((\a+n+1)F_{1,n}(k_0-1,1)-\frac{1}{\a+n+1}$$
  $$-F_{1,1+n}(k_0-1,1)(\a+n))-
  \frac{4}{(\a+n+4)^2}$$
  $$\hskip -4cm =\frac{\a+n+2}{\a+n+4}((\a+n+1)F_{1,n}(k_0-1,1)-\frac{1}{\a+n+1}$$
$$
-(\frac{\a+n+1}{\a+n+2}F_{1,n}(k_0-1,1)-\frac{1}{(\a+n+2)^2})(\a+n))-
  \frac{4}{(\a+n+4)^2},$$
and
  $$F_{2,n+2}(k_0-1,1)=(\a+n+3)F_{1,n+2}(k_0-1,1)-\frac{1}{\a+n+3}-F_{1,3+n}(k_0-1,1)(\a+n+2)$$
  $$\hskip -5cm=(\a+n+3)F_{1,n+2}(k_0-1,1)-\frac{1}{\a+n+3}$$
  $$-(\a+n+2)(\frac{\a+n+3}{\a+n+4}F_{1,n+2}(k_0-1,1)-\frac{1}{(\a+n+4)^2})$$
  $$=(\a+n+3)\frac{2}{\a+n+4}F_{1,n+2}(k_0-1,1)-\frac{1}{\a+n+3}+\frac{\a+n+2}{(\a+n+4)^2}$$
  $$\hskip -2cm=\frac{2(\a+n+3)}{\a+n+4}(\frac{\a+n+2}{\a+n+3}(\frac{\a+n+1}{\a+n+2}F_{1,n}(k_0-1,1)
  $$
  $$-\frac{1}{(\a+n+2)^2})-\frac{1}{(\a+n+3)^2})-\frac{1}{\a+n+3}+\frac{\a+n+2}{(\a+n+4)^2}.$$
Equating the above two expressions for $F_{2,n+2}(k_0-1,1)$ and
simplifying yield
$${\frac {-6}{ \left( \a+n+3 \right)  \left( \a+n+4 \right) ^{2}
 \left( \a+n+1 \right) }}=0,\,\,\forall\,\,n\in\Z,$$
a contradiction. Then $F_{1,0}({k_0},1)=0$, furthermore
 $F_{j,n}({k_0},1)=0 \,\,\forall\,\,j, n\in \Z.$

 \vskip 5pt
{\it Case 1.2.2:
$A_{\pm1,n}(k_0-1,k_0)=0,A_{\pm2,n}(k_0-1,k_0)=\pm
\frac{1}{(\a+n\pm 1)(\a+n\pm2)}.$}

 Again using (4.9) and (4.10), from the $(k_0-1,1)$-entry of (4.4), we
 obtain
$$(\a+n+j)F_{j,n}(k_0-1,1)+A_{i,j+n}(k_0-1,k_0)\frac{j}{\a+n+j}-F_{j,i+n}(k_0-1,1)(\a+n)$$
$$=jF_{i+j,n}(k_0-1).\eqno(4.14)$$ Letting $j=1$ and $i=-1$, we obtain that
$$0=(\a+n+1)F_{1,n}(k_0-1,1)+A_{-1,1+n}(k_0-1,k_0)\frac{1}{\a+n+1}-F_{1,n-1}(k_0-1,1)(\a+n),$$
i.e.
$$F_{1,n}(k_0-1,1)=F_{1,n-1}(k_0-1,1)\frac{\a+n}{\a+n+1}, \,\,\forall\,\,j, n\in \Z.\eqno(4.15)$$
Letting $i=j=1$, and $j=-i=2$ in (4.14) respectively, we obtain
that
$$F_{2,n}(k_0-1,1)=(\a+n+1)F_{1,n}(k_0-1,1)-F_{1,n+1}(k_0-1,1)(\a+n),$$
$$(\a+n+2)F_{2,n}(k_0-1,1)=\frac{1}{(\a+n)(\a+n+1)}\frac{2}{\a+n+2}+F_{2,n-2}(k_0-1,1)(\a+n).$$
Using the above two formulas and (4.15) we deduce that
$$\hskip -8cm F_{2,n+2}(k_0-1,1)$$
$$=\frac{2}{(\a+n+2)(\a+n+3)(\a+n+4)^2}+F_{2,n}(k_0-1,1)\frac{\a+n+2}{\a+n+4}$$
$$\hskip -6cm=\frac{2}{(\a+n+2)(\a+n+3)(\a+n+4)^2}$$ $$+((\a+n+1)F_{1,n}(k_0-1,1)
-(F_{1,n}(k_0-1,1)\frac{\a+n+1}{\a+n+2})(\a+n))\frac{\a+n+2}{\a+n+4},$$
and
 $$F_{2,n+2}(k_0-1,1)=(\a+n+3)F_{1,n+2}(k_0-1,1)-F_{1,n+3}(k_0-1,1)(\a+n+2)$$
 $$=(\a+n+1)F_{1,n}(k_0-1,1)-F_{1,n}(k_0-1,1)(\a+n+1)(\a+n+2)/(\a+n+4).$$
Equating the above two expressions for $F_{2,n+2}(k_0-1,1)$ and
simplifying yield that $0=\frac{2}{(\a+n+2)(\a+n+3)(\a+n+4)^2}$
for all $n$, a contradiction. Then $F_{1,0}({k_0},1)=0$,
furthermore $F_{j,n}({k_0},1)=0 \,\,\forall\,\,j, n\in \Z.$

\medskip
{\it Case 2:  $k_0=2$.}

{\it Case 2.1:  $V$ is decomposable over $\Vir$. }

\vskip 5pt
\medskip From the established Case 1 we may assume that $V=W^{(2)}$.
Note that $A_{i,n}(1,2)$ $=0$. The $(1,1)$-entry of (4.4) gives
$$(\a+j+n)F_{j,n}(1,1)-F_{j,i+n}(1,1)(\a+n)=jF_{i+j,n}(1,1)+\delta_{i,-j}(i^2+i)c_{DI}.\eqno(4.16)$$
 Letting $j=1$, we obtain
 $$F_{i+1,n}(1,1)=(\a+1+n)F_{1,n}(1,1)-F_{1,i+n}(1,1)(\a+n)-\delta_{i,-1}(i^2+i)c_{DI}.\eqno(4.17)$$
Then letting $i=-1$, we have
$$F=(\a+1+n)F_{1,n}(1,1)-F_{1,n-1}(1,1)(\a+n).$$

If $F=0$, then we have
$(\a+1+n)F_{1,n}(1,1)=F_{1,n-1}(1,1)(\a+n)$, i.e.,
$$F_{1,n}(1,1)=(\a+1)F_{1,0}(1,1)/(\a+1+n).$$ Applying this to
(4.17) we obtain that
$$F_{j,n}(1,1)=\frac{(\a+1)j}{\a+n+1}F_{1,0}(1,1),\eqno(4.18)$$ even for $j=2$ (Using (4.16)).

Since $W^{(2)}$ is decomposable, by symmetry of (4.9) and (4.18)
we have
$$F_{j,n}(1,2)=\frac{(\a+1)j}{\a+n+j}F_{1,0}(1,2),\,\,\,\forall\,\,\,j, n\in\Z.\eqno(4.19)$$
$$F_{j,n}(2,2)=\frac{(\a+1)j}{\a+n+1}F_{1,0}(2,2).\eqno(4.20)$$
Thus
$$F_{j,n}=\frac{(\a+1)j}{\a+n+j}\pmatrix F_{1,0}(1,1)& F_{1,0}(1,2)\\ F_{1,0}(2,1) &
F_{1,0}(2,2)\endpmatrix ,\,\,\,\forall\,\,\,j,
n\in\Z.\eqno(4.21)$$ Since the constant matrix has an eigenvector,
there exist $a,b\in\bC$ such that $W=\sum\bC(av_n^1+bv_n^2)$ is a
nontrivial ${\Bbb L}$-submodule, a contradiction.

If $F\neq 0$, we have
$F_{j,n}(1,1)-F=j\frac{\a+1}{\a+n+1}(F_{1,0}(1,1)-F)$. Similarly
we have
$$F_{j,n}=\frac{(\a+1)j}{\a+n+j}\pmatrix F_{1,0}(1,1)-F& F_{1,0}(1,2)\\ F_{1,0}(2,1) &
F_{1,0}(2,2)-F\endpmatrix +\pmatrix F& 0\\
0 & F\endpmatrix,\,\,\,\forall\,\,\,j, n\in\Z,$$ again we have a
contradiction.

\medskip
{\it Case 2.2:  $W^{(2)}$ is indecomposable  over $\Vir$.}
\medskip The argument  is exactly the same as in Case 1.2.2. We do
not repeat it.

So far we have proved that $F_{j,n}(k,1)=0$ for all $n ,j\in\Z$
and $k\ne1$. Thus $W^{(1)}$ is an $\Bbb L$-submodule which must be
$V$. Combining with Lemma 3.3, we have proved this lemma.
 \hfill $ $\qed
\enddemo

Denote by $(W^{(i)}/W^{(i-1)})'$ the unique nontrivial
sub-quotient $\Vir$-module of\break $W^{(i)}/W^{(i-1)}$. For any
$x,y \in \bC$, define $x\nprec y$ if $y-x \notin \N $.

\proclaim{Lemma 4.3} The module $V$ carries a filtration
$\{W^{(1)},W^{(2)}, \cdots ,W^{(p)}\}$ with 
$(W^{(i)}/W^{(i-1)})'\simeq V'(\a,\b_i)$ as $\Vir$-modules, where
$ \b_i \nprec \b_j$ for all $i<j$.
\endproclaim
\demo{Proof}: We start with the filtration at the beginning of
this section. The statement is true if $p=1$, or $2$ (use Theorem
2.2 and $V'(\a,0)\simeq V'(\a,1)$). Now we consider $p>2$.

Suppose that we do not have $ \b_{i} \nprec \b_{i+1}$ for some
$i$, say $ \b_{p-1} \nprec \b_{p}$, i.e., $\b_{p}-\b_{p-1}\in \N$.
 Consider $V/V^{(p-2)}$ ($p=2$
for this module). Then we can have a submodule $ X^{(p-1)}\supset
V^{(p-2)}$ such that $(W^{(p)}/X^{(p-1)})'\simeq V'(\a,\b_{p-1})$
and $(X^{(p-1)}/W^{(p-2)})'\simeq V'(\a,\b_{p})$.

By repeating this procedure several times  if necessary, then we
obtain the filtration required. \hfill $ $\qed
\enddemo
Now we are ready to classify all irreducible uniformly bounded
weight  modules over $\Bbb L$.

 \proclaim{Theorem 4.4} If V is a nontrivial
irreducible uniformly bounded weight  module over $\Bbb L$, then V
is isomorphic to $V'(\a,\b;F)$ for some $\a,\b,F \in \bC$.
\endproclaim

\demo{Proof} By Lemma 4.3, the module V carries a filtration
$\{W^{(1)},W^{(2)}, \cdots ,W^{(p)}\}$ with
$(W^{(i)}/W^{(i-1)})'\simeq V'(\a,\b_i)$ such that $$ \b_i \nprec
\b_j \,\,\, {\text{ for}}\,\,\,   {\text{ all}}\,\,\,  i<j. \tag
4.22$$

Suppose  that $\b_i=0,1$ for all $i$. Since $V'(\a,0)\simeq
V'(\a,1)$, then by  Lemmas 4.1 and 4.2 we need  consider only the
case $\a =0$ and $c_{DI}\neq F$. In this case, consider suitable
Vir[$e$] instead of $\Vir[0]$, we have $\a \notin \Z$ which
applies to Lemma 4.2. Thus the only remaining  case is the case
$\b_i \neq 0,1 $ for some $i$. By replacing $V$ with its
contragredient module if necessary, from (4.22) we may assume that
$\b_1 \neq 0,1.$

\medskip
Suppose $F_{j,n}(k,1)=0$ for all $j,n\in\Z$, $k>k_0$, where
$k_0>1$ is a fixed integer. We need to show that
$F_{j,n}(k_0,1)=0$ for all $j,n\in\Z$.

\medskip
{\bf Claim 1.} {\it $F_{1,n}(k_0,1)=0$  except for finitely many
$n\in\Z$.}

\medskip

{\it Case 1: $\a \notin \Z$}.

 In this case all the restrictions for (4.2)-(4.4) disappear. Then the $(k_0,1)$-entry of  (4.4) gives
 $$(\a+n+j+i\b_{k_0})F_{j,n}(k_0,1)-F_{j,i+n}(k_0,1)(\a+n+i\b_1)=jF_{i+j,n}(k_0,1).\eqno(4.23)$$
 Letting $j=1$ we obtain
$$F_{i+1,n}(k_0,1)=(\a+n+1+i\b_{k_0})F_{1,n}(k_0,1)-(\a+n+i\b_1)F_{1,i+n}(k_0,1).\eqno(4.24)$$
Taking $i=-1$   we have
$$
(\a-\b_{k_0}+n+1)F_{1,n}({k_0},1)=(\a-\b_1+n)F_{1,n-1}({k_0},1).\eqno(4.25)$$
 Letting $i=1$ in (4.24), we have
 $F_{2,n}({k_0},1)=(\a+\b_{k_0}+n+1)F_{1,n}({k_0},1)-(\a+\b_1+n)F_{1,n+1}({k_0},1).$
 So
$$(\a-\b_{k_0}+n+2)F_{2,n}({k_0},1)=(\a-\b_{k_0}+n+2)(\a+\b_{k_0}+n+1)F_{1,n}({k_0},1)$$
$$-(\a+\b_1+n)(\a-\b_1+n+1)F_{1,n}({k_0},1)$$
$$=(2\a-\b_{k_0}^2+\b_{k_0}+2n+\b_1^2-\b_1+2)F_{1,n}(k_0,1). \eqno(4.26)$$
Applying (4.25) to (4.26) we deduce that
$$(\a-\b_{k_0}+n+2)(\a-\b_{k_0}+n+3)(\a-\b_{k_0}+n+4)F_{2,n+2}({k_0},1)$$
$$=(\a-\b_1+n+1)(\a-\b_1+n+2)(2\a-\b_{k_0}^2+\b_{k_0}+2n+4+\b_1^2-\b_1+2)F_{1,n}(k_0,1).\eqno(4.27)$$
By using  (4.23) with $i=-2$, $ j=2$, and $n$ being replaced by
$n+2$, we deduce that
$$(\a+n+4-2\b_{k_0})(\a-\b_{k_0}+n+2)F_{2,2+n}(k_0,1)$$
$$=(\a+n+2-2\b_1)(\a-\b_{k_0}+n+2)F_{2,n}(k_0,1)$$
$$=(\a+n+2-2\b_1)(2\a-\b_{k_0}^2+\b_{k_0}+2n+\b_1^2-\b_1+2)F_{1,n}(k_0,1).\eqno(4.28)$$
Combining (4.27) and (4.28), we deduce that a product is $0$, one
factor of the product is $F_{1,n}(k_0,1)$, and the other is
$$(\a+n+4-2\b_{k_0})(\a-\b_1+n+1)(\a-\b_1+n+2)(2\a-\b_{k_0}^2+\b_{k_0}+2n+4+\b_1^2-\b_1+2)$$
$$-
(\a-\b_{k_0}+n+3)(\a-\b_{k_0}+n+4)(\a+n+2-2\b_1)(2\a-\b_{k_0}^2+\b_{k_0}+2n+\b_1^2-\b_1+2),$$
for all $n\in\Z$.  From simplifying this  we obtain that
$$({{\b_{k_0}}}-{\b_1}-1)({{\b_{k_0}}}-{\b_1}-2)[\left( { \b_{k_0}}^{2}-3\,{ \b_{k_0}}+2\,{ \b_{k_0}}\,{ \b_1}-{ \b_1}+{{
 \b_1}}^{2} \right) (n+\a)$$
 $$-2\, \left( {\b_1}-1 \right)  \left( {\b_{k_0}}
\,{\b_1}-2\,{\b_1}+{{\b_{k_0}}}^{2}-3\,{\b_{k_0}}
\right)]F_{1,n}(k_0,1)=0.$$ Using (4.22) and the fact that
$\b_1\ne0, 1$, we deduce that
$$F_{1,n}(k_0,1)=0\eqno(4.29)$$ for all but one possible $n=n_0$, and from
(4.25) we know that $\a-\b_1+n_0+1=0$.
Claim 1 follows in this case.

\medskip
 {\it Case 2:  $\a=0$}.
\medskip

Since $W^{(1)}\simeq V(0,\b_1)$ and $\b_1\ne0$ or $1$, then
$\dim(W^{(1)})_0=1$ and we can have $v_0^1$. In this case, the
restrictions for (4.4) become $(n+j)(n+i+j)\ne0.$ The restrictions
for (4.23)-(4.28) are $(n+j)(n+i+j)\ne0$, $(n+1)(n+i+1)\ne0$,
$n(n+1)\ne0$, $(n+1)(n+2)\ne0$, $(n+1)(n+2)(n+3)(n+4)\ne0$ and
$(n+1)(n+2)(n+3)(n+4)\ne0$, respectively. Thus we  have (4.29)
with exceptions $n=-1,-2,-3,-4$ and possibly one more exception
$n=n_0$ (which comes from the computation of getting (4.29)).
Claim 1 follows.

\medskip
{\bf Claim 2.} {\it There exists some $i_0 \in \Z, i_0 \neq -1,0$,
such that $I(i_0+1)W^{(1)} \subseteq W^{(k_0-1)}$. }\medskip

  For any $j\in\Z\setminus\{0\}$, set
$$S_j=\{ n\in \Z | I(j)v_{n}^1 \nsubseteq W^{(k_0-1)}\}.$$
By Claim 1 we know that $|S_1|< +\infty$. Choose $i_0$ satisfying

    (a) $i_0 \neq -1,0$;

    (b) $i_0> \max\{|x-y| | x, y \in S_1\}+1$;

    (c) $-\a+(-1+\b_1)i_0 \notin S_1, i_0\b_1+\b_1-\a-1 \notin S_1$ (because $\b_1 \neq
    0,1$).
 \vskip 5pt
 It is clear that
    $$(i_0k+S_1)\cap S_1=\emptyset,\,\,\,\text{for }\,\,\, k\ne0.$$
 From
    $$(\a+n-i_0+i_0\b)I(1)v_{n}^1= I(1)D(i_0)v_{n-i_0}^1=-I(i_0+1)v_{n-i_0}^1+D(i_0)I(1)v_{n-i_0}^1,$$
    we have
$$I(i_0+1)v_{n-i_0}^1 \in W^{(k_0-1)}\,\,\,\text{if}\,\,\,n \notin S_1\cup( i_0+S_1).\eqno(4.30)$$
Since $D(-i_0)I(i_0+1)v_{n-i_0}^1 \in W^{(k_0-1)}$ for all $n
\notin S_1\cup(     i_0+S_1)$, i.e.,
      $$(i_0+1)I(1)v_{n-i_0}^1+I(i_0+1)D(-i_0)v_{n-i_0}^1$$
      $$=(i_0+1)I(1)v_{n-i_0}^1+(\a+n-i_0-i_0\b)I(i_0+1)v_{n-2i_0}^1 \in W^{(k_0-1)},\,\,\,\forall \,\,\,n \notin S_1\cup(
      i_0+S_1),$$
we deduce that $$(\a+n-i_0-i_0\b)I(i_0+1)v_{n-2i_0}^1 \in
      W^{(k_0-1)},\,\,\,\forall \,\,\,n \notin S_1\cup(i_0+S_1).\eqno(4.31)$$
      For any $n \in 2i_0+S_1$, from (b) we have $n \notin S_1$ and
      $n \notin i_0+S_1,$  and from (c), we have $\a+n-i_0-i_0\b_1\neq 0$.
      Applying this to (4.31) we obtain that
      $$I(i_0+1)v_{\a+n}^1 \in
      W^{(k_0-1)},\,\,\,\forall \,\,\,n \in S_1.$$
      Together with (4.30), we have
   $$I(i_0+1)v_{n-i_0}^1 \in W^{(k_0-1)},\,\,\,\forall \,\,\,n\notin
   S_1.\eqno(4.32)$$
From $[D(-i_0-1),I(i_0+1)]=(i_0+1)I(0)+((-i_0-1)^2-i_0-1)c_{DI}$,
      we have $D(-i_0-1)I(i_0+1)v_{n-i_0}^1\equiv
      I(i_0+1)D(-i_0-1)v_{n-i_0}^1 \mod W^{(k_0-1)}.$
      Hence
     $$(\a+n-i_0-(i_0+1)\b_1)I(i_0+1)v_{n-2i_0-1}^1 \in
      W^{(k_0-1)},\,\,\,\forall \,\,\,n\notin S_1.\eqno(4.33)$$
     For $n-i_0-1\in S_1$ we know that $n\notin S_1$, and  by (c) we also have $(\a+n-i_0-(i_0+1)\b_1)\ne0$. Thus
     $$I(i_0+1)v_{n-i_0}^1 \in
      W^{(k_0-1)},\,\,\,\forall \,\,\,n\in S_1.$$
      Therefore $S_{i_0+1}=\emptyset.$
      So we have proved this Claim.

 Noting that $\{ x_i, I(i_0)|i\in\Z\}$ generate  ${\Bbb L}$, we
 know that ${\Bbb L}W^{(1)}\subset W^{(k_0-1)}$. By induction on
 $k_0$ we see that $W^{(1)}$ is an $\Bbb L$-submodule. From
Theorem 3.3 we complete the proof of the theorem. \hfill \qed
\enddemo

\subhead 5. \ \  Classification of irreducible weight modules over
$\Bbb L$ with finite-dimensional weight spaces
\endsubhead
\medskip
\proclaim{Theorem 5.1} Let  V be an irreducible weight module over
$\Bbb L$ with all weight spaces finite-dimensional. If V is not
uniformly bounded, then V is either a highest weight module or a
lowest weight module.

\demo{Proof} Consider  $V$ as a $\Vir[0]$-module. Let $W$ be the
smallest $\Vir[0]$-submodule of $V$  such that $V/W$ is a trivial
$\Vir[0]$-submodule. Then $W$ contains no trivial quotient module.

Let $W'$ be the maximal trivial $\Vir[0]$-submodule of $W$. Then
$W/W'$ contains no trivial submodule.

Since $\dim W'+\dim V/W$ is finite, then $\Vir[0]$-module $W/W'$
is not uniformly bounded. Now $W/W'$ satisfies the conditions in
Theorem 2.4. By using Theorem 2.4 to $W/W'$, we have some
nontrivial upper bounded (or lower bounded) $\Vir[0]$-submodule of
$W/W'$, say, $W''/W'$ is a nontrivial upper bounded
$\Vir[0]$-submodule of $W/W'$ (i.e., the weight set of $W''/W'$
has an upper bound).

Denote $W''=M^+$, We know that $M^+$ is not uniformly bounded.

For any $j\in\N$, define $M^+(j)=\{v\in M^+|I(i)
v=0\,\,\forall\,\, i\ge j\}$, and let
$$M=\cup_{i\in\N} M^+(j).$$
It is easy to check that $x_i M^+(j)\subset M^+(j+|i|)$, i.e., $M$
is an   $\Vir[0]$-submodule of $M^{+}$. Suppose that
$\supp(V)\subset \a+\Z$ for some $\a\in\bC$.

\medskip
{\bf Claim. $M\ne0$.}

Fix $\lambda_0 \in \a+\Z.$ Since $ M^+$ is not uniformly bounded,
we have some $0\neq \lambda_1 \in \supp(M^+)$ with
$\lambda_1<\lambda_0$ and  $\dim (M^+)_{\lambda_1}
>\dim V_{\lambda_0}$. Hence
$I(\lambda_0-\lambda_1):(M^+)_{\lambda_1}\rightarrow
V_{\lambda_0}$ is not injective. Say $v=v_{\lambda_1} \in
M^+\setminus\{0\}$ with $I(i_0)v=0$ where
$i_0=\lambda_0-\lambda_1>0$. Since $v\in M^+$, there exists
$j_0>0$ such that  $x_{j}v=0$ for $j\ge j_0$ Then we have
$I(i)v=0$ for all $i\ge i_0+j_0$. Thus $v\in M$. Claim follows.
\medskip
 Let $\Lambda$ be the maximal weight of $M$, and
$v_{\Lambda}$ is one of the corresponding weight vectors. By the
definition of $M$,  there exists a nonnegative integer $i_0$ such
that $I(i)v_{\Lambda}= 0$ for $i>i_ 0$, and
$I(i_0)v_{\Lambda}\ne0$ if $i_0>0$. If $i_0=0$, then $v_{\Lambda}$
is a highest weight vector of the ${\Bbb L}$-module $V$, and we
are done. So we assume that $i_0>0$. From
$x_jI(i_0)v_{\Lambda}=i_0I(i_0+j)v_{\Lambda}+I(i_0)x_jv_{\Lambda}=I(i_0)x_jv_{\Lambda}=0$
for all $j>0$, we know that $I(i_0)v_{\Lambda}\ne0 $ is a  highest
weight vector over $\Vir[0]$. So $I(i_0)v_{\Lambda}\in M^+$. From
$I(i)I(i_0)v_{\Lambda}=I(i_0)I(i)v_{\Lambda}=0$ for all $i>i_0$,
we know that $I(i_0)v_{\Lambda}\in M$, contradicting the choice of
$\Lambda$ and $i_0>0$. So we have proved this Theorem. \hfill $
$\qed
\enddemo

 Combining Theorems 4.4 and 5.1 we obtain the main result of this paper:

\proclaim{Theorem 5.2} If V is a nontrivial irreducible weight
module over $\Bbb L$ with finite dimensional weight spaces, then V
is isomorphic to $V'(\a,\b;F)$ for some $\a,\b, F \in \bC$, or  a
highest or lowest weight module.
\endproclaim

\vskip.3cm \Refs\nofrills{\bf REFERENCES}
\bigskip
\parindent=0.45in

\leftitem{[A]} E. Arbarello, C. De Concini, V. G. Kac and C.
Procesi, Moduli spaces of curves and representation theory, {\it
Comm. Math. Phys.}, (1)117, 1-36(1988).

\leftitem{[B]} Y. Billig,  Representations of the twisted
Heisenberg-Virasoro algebra at level zero, {\it Canad. Math.
Bull.}, 46(2003), no.4, 529-537.

\leftitem{[KR]} V. G. Kac and K. A. Raina, ``Bombay lectures on
highest weight representations of infinite dimensional Lie
algebras,'' World Sci., Singapore, 1987.

\leftitem{[M]} O. Mathieu, Classification of Harish-Chandra
modules over the Virasoro algebra, {\it Invent. Math.} {\bf
107}(1992), 225-234.

\leftitem{[Ma1]} V. Mazorchuk, Classification of simple
Harish-Chandra modules over ${\bQ}$-Virasoro algebra, {\it Math.
Nachr.} {\bf 209}(2000), 171-177.

\leftitem{[Ma2]} V. Mazorchuk,  Verma modules over generalized
Witt algebras, {\it Compositio Math.} (1) {\bf 115}(1999), 21-35.

\leftitem{[MP1]} C. Martin and A. Piard, Nonbounded Indecomposable
Admissible Modules over the Virasoro Algebra, {\it Lett. Math.
Phys.} 23, 319-324(1991).

\leftitem{[MP2]} C. Martin and A. Piard, Indecomposable Modules
Over the Virasoro Lie algebra and a Conjecture of Kac, {\it
Commun. Math. Phys.} 137, 109-132(1991).

\leftitem{[MP3]} C. Martin and A. Piard, Classificaion of the
Indecomposable Bounded Admissible Modules over the Virasoro Lie
algebra with Weightspaces of Dimension not Exceeding Two, {\it
Commun. Math. Phys.} 150, 465-493(1992).


\leftitem{[SZ]} Y. Su and K. Zhao, Generalized Virasoro and
super-Virasoro algebras and modules of intermediate series, {\it
J. Algebra}, 252(2002), no.1, 1-19.

\leftitem{[Zh]} K. Zhao, Weight modules over generalized Witt
algebras with $1$-dimensional weight spaces, {\it Forum. Math.},
 Vol.16, No.5, 725-748(2004).
\endRefs
\vfill
\enddocument
\end